\documentclass [12pt]{article}
\title {Definability and continuity of the SU-rank in unidimensional supersimple theories}
\author {Ziv Shami}
\newtheorem {theorem}{Theorem}[section]
\newtheorem {lemma}[theorem]{Lemma}
\newtheorem {definition}[theorem]{Definition}

\newtheorem {fact}[theorem]{Fact}

\newtheorem {corollary}[theorem]{Corollary}

\newtheorem {remark}[theorem]{Remark}

\newtheorem {notation}[theorem]{Notation}
\newtheorem {claim}[theorem]{Claim}
\newtheorem {subclaim}[theorem]{Subclaim}

\def\proof {\noindent \textbf{Proof:} }
\def\proofof #1 {\noindent \textbf{Proof of {#1}:} }
\def\qed {$\ \ \ \ \Box$}

\def\qedend {$\ \ \ \ \ \ \ \ \ \ \ \ \ \ \ \ \ \ \ \ \ \ \ \ \ \ \ \ \ \ \ \ \ \ \ \ \ \ \ \ \ \ \ \ \ \ \ \ \ \ \
 \ \ \ \ \ \ \ \ \ \ \ \ \ \ \ \ \ \ \ \ \ \ \ \ \ \ \ \ \ \ \ \ \ \ \ \ \ \ \ \ \ \ \ \Box$\\}

\def\Ceq  {\CC^{eq}}

\newsavebox{\indbin}
\savebox{\indbin}{\begin{picture}(0,0)
\newlength{\gnu}
\settowidth{\gnu}{$\smile$} \setlength{\unitlength}{.5\gnu} \put(-1,-.65){$\smile$}
\put(-.25,.1){$|$}
\end{picture}}
\newcommand{\nonfork}[3]
{\mbox{$\begin{array}{ccc} \mbox{$#1$} & \usebox{\indbin} & \mbox{$#2$} \\
        & \mbox{$#3$} &
\end{array}$}}

\newcommand{\fork}[3]
{\mbox{$\begin{array}{ccc} \mbox{$#1$} & \!\mbox{$\!\!\not\!\:\usebox{\indbin}$} & \mbox{$#2$} \\
        & \mbox{$#3$} &
\end{array}$}}

\newsavebox{\sindbin}
\savebox{\sindbin}{\begin{picture}(0,0)
\newlength{\sgnu}
\settowidth{\sgnu}{$\smile$} \setlength{\unitlength}{.5\sgnu} \put(-1,-.65){$\smile$}
\put(-.25,.1){$|s$}
\end{picture}}

\newsavebox{\starindbin}
\savebox{\starindbin}{\begin{picture}(0,0)
\newlength{\stargnu}
\settowidth{\stargnu}{$\smile$} \setlength{\unitlength}{.5\stargnu} \put(-1,-.65){$\smile$}
\put(-.25,.1){$|*$}
\end{picture}}

\newsavebox{\qindbin}
\savebox{\qindbin}{\begin{picture}(0,0)
\newlength{\qgnu}
\settowidth{\qgnu}{$\smile$} \setlength{\unitlength}{.5\qgnu} \put(-1,-.65){$\smile$}
\put(-.25,.1){$|_{qf}$}
\end{picture}}

\def\card #1 {{\vert #1 \vert}}

\def\CC {{\cal C}}

\def\PP {{\cal P}}
\def\UU {{\cal U}}

\usepackage{amsfonts}

\begin{document}
\maketitle

\begin{abstract}
We prove, in particular, that in a supersimple unidimensional theory the $SU$-rank is continuous
and the $D$-rank is definable.
\end{abstract}

\section{Introduction}
Analysis of "Forking sets", i.e. invariant sets built by applying first-order operations on the
predicates $R(a,b)$ if "$\phi(x,a)\mbox{\ forks\ over}\ b$", has powerful implications on the
structure of definable sets in simple theories, e.g. Lstp=stp for low theories, elimination of
hyperimaginaries in supersimple theories, supersimplicity of countable hypersimple unidimensional
theories and of a large class (non s-essentially 1-based theories) of uncountable hypersimple
unidimensional theories and more (a theory is \em hypersimple \em if it is simple and eliminates
hyperimaginaries).

It is well known that a unidimensional stable theory is superstable. Stable unidimensional theories
share some nice definability properties: the $SU$-rank is continuous and Shelah's $D$-rank is
definable. In this paper, we generalize and strengthen these properties to simple theories. We
prove this by an argument that relies on this kind of analysis. We will assume basic knowledge of
simple theories as in [K1],[KP],[HKP]. A good text book on simple theories that covers much more is
[W]. The notations are standard, and throughout the paper we work in a highly saturated, highly
strongly-homogeneous model $\CC$ of a complete first-order theory $T$ in a language $L$.

\section{Preliminaries}
In this section $T$ is assumed to be simple. We quote several known facts that we will apply.

\subsection{Lowness}
We will say that the formula $\phi(x,y)\in L$ is \em low in $x$ \em if there exists $k<\omega$ such
that for every $\emptyset$-indiscernible sequence $(b_i \vert i<\omega)$, the set $\{\phi(x,b_i)
\vert i<\omega\}$ is inconsistent iff every subset of it of size $k$ is inconsistent. Note that
every stable formula $\phi(x,y)$ is low in both $x$ and $y$.

\begin{remark}\label{low type def}\em
Note that if $\phi(x,y)\in L$ is low in $x$ then the relation $F_\phi$ defined by $F_\phi(b,A)$ iff
$\phi(x,b)$ forks over $A$ is type-definable.
\end{remark}

\subsection{Almost internality and analyzability}

Let $\PP$ be an $A$-invariant family of partial types. We say that \em $p \in S(A)$ is (almost-)
$\PP$-internal \em if there exists a realization $a$ of $p$ and there exists $b$ with
$\nonfork{a}{b}{A}$ such that for some tuple $c$ of realizations of partial types in $\PP$ over
$Ab$ we have $a \in dcl(b,c)$ (respectively, $a \in acl(b,c)$). We say that \em $p$ is analyzable
in $\PP$ \em if there exists a sequence $I=\langle a_i \vert i\leq \alpha \rangle$ in $\CC^{eq}$
such that $a_\alpha\models p$ and $tp(a_i/\{a_j \vert j<i\}\cup
A)$ is almost $\PP$-internal for every $i\leq \alpha$.\\

\noindent First, the following fact is straightforward.

\begin{fact}\label{some internal}
1) Assume $tp(a_i/A)$ are (almost) $\PP$-internal for $i<\alpha$. Then $tp(\langle a_i \vert
i<\alpha \rangle/A)$ is (almost) $\PP$-internal. Thus, if $tp(a_i/A)$ are analyzable in $\PP$ for
$i<\alpha$. Then $tp(\langle a_i \vert i<\alpha \rangle/A)$ is analyzable in $\PP$.

\noindent 2) If $tp(a/A)$ almost $\PP$-internal, so is $tp(a/B)$ for any set $B\supseteq A$.
\end{fact}

\noindent The following Fact [S0, Theorem 2.2] will useful (a quite similar result has been proved
independently in [W, Proposition 3.4.9]).

\begin {fact}\label{fact_a_int}
Let $T$ be simple. Let $p(x)\in S(A)$ be a hyperimaginary amalgamation base and let $\UU$ be an
$A$-invariant set of hyperimaginaries. Suppose $p$ is almost $\UU$-internal. Then for every Morley
sequence $\bar a$ of length $\omega$ of $p$ there is a type-definable one-to-bounded relation
$S(x,\bar y,\bar a)$ (i.e. for every $\bar y$ there are boundedly many x-s for which $S(x,\bar
y,\bar a)$ holds) which covers p by $\UU$. If $p$ and $\UU$ are real then $S$ can be chosen to be
definable.
\end {fact}

\subsection{The extension property}

We recall some natural extensions of notions from [BPV]. By a pair $(M,P^M)$ of $T$ we mean an
$L_P=L\cup \{P\}$-structure, where $M$ is a model of $T$ and $P$ is a new predicate symbol whose
interpretation is an elementary submodel of $M$.  For the rest of this subsection, by a $\vert T
\vert$-\em small type \em we mean a complete hyperimaginary type in $\leq \vert T\vert$ variables
over a hyperimaginary of length $\leq \vert T\vert$ (i.e. a sequence of length $\leq \vert T\vert$
modulo a $\emptyset$-type-definable equivalence relation).

\begin{definition}\label {def1}\em
Let $\PP_0,\PP_1$ be $\emptyset$-invariant families of $\vert T\vert$-small types.

\noindent 1) We say that a pair $(M,P^M)$ satisfies the extension property for $\PP_0$ if for every
$L$-type $p\in S(A)$, where $A$ is a hyperimaginary with $A\in dcl(M)$ and $p\in \PP_0$, there is
$a\in p^M$ such that $\nonfork{a}{P^M}{A}$.

\noindent 2) Let $$T_{Ext,\PP_0}=\bigcap \{Th_{L_P}(M,P^M) \vert\ \mbox{the\ pair\ $(M,P^M)$
satisfies\ the\ extension\ property\  w.r.t.\  } \PP_0\ \}.$$

\noindent 3) We say that $\PP_0$ dominates $\PP_1$ w.r.t. the extension property if $(M,P^M)$
satisfies the extension property for $\PP_1$ for every $\vert T\vert^+$-saturated pair
$(M,P^M)\models T_{Ext,\PP_0}$. In this case we write $\PP_0\unrhd_{_{Ext}} \PP_1$.

\noindent 4) We say that the extension property is first-order for $\PP_0$ if
$\PP_0\unrhd_{_{Ext}}\PP_0$ (i.e. every $\vert T\vert^+$-saturated model of $T_{Ext,\PP_0}$
satisfies the extension property for $\PP_0$). We say that the extension property is first-order if
the extension property is first-order for the family of all $\vert T\vert$-small types
(equivalently, for the family of all real types over sets of size $\leq \vert T\vert$).
\end{definition}

\begin {fact} \label {foext} [BPV, Proposition\ 4.5] \em
The extension property is first-order in $T$ iff for every formulas $\phi(x,y),\psi(y,z)\in L$ the
relation $Q_{\phi,\psi}$ defined by: $$Q_{\phi,\psi}(a)\mbox{\ iff}\ \phi(x,b)\mbox{ doesn't\ fork\
over}\ a\ \mbox{for\ every}\ b\models \psi(y,a)$$ is type-definable (here $a$ can be an infinite
tuple from $\CC$ whose sorts are fixed).
\end {fact}

\begin{fact} \label{ext family}[S1, Lemma 3.7]
Let $\PP_0$ be an $\emptyset$-invariant family of $\vert T\vert$-small types. Assume $\PP_0$ is
extension-closed and that the extension property is first-order for $\PP_0$. Let $\PP^*$ be the
maximal class of $\vert T\vert$-small types such that $\PP_0\unrhd_{_{Ext}}\PP^*$. Then
$\PP^*\supseteq An(\PP_0)$, where $An(\PP_0)$ denotes the class of all $\vert T\vert$-small types
analyzable in $\PP_0$ by hyperimaginaries.
\end{fact}

We will need the following consequence.

\begin{corollary}\label{epfo_su1_formula}
Assume $\theta=\theta(x)\in L$, $SU(\theta)=1$ and $tp(a)$ is analyzable in $\theta$ for every
$a\in\CC$. Then the extension property is first-order in $T$.
\end{corollary}

\proof First, note the following slightly stronger version of Hrushovski's Lemma [H, Lemma 4.3]:

\begin{claim}\label{H_claim}
Let $\theta(x)\in L$, $SU(\theta(x))=1$ and let $\chi(x,y)\in L$ be such that $\chi(x,y)\vdash
\theta(x)$. Then there exists $N<\omega$ such that for all $a\in \CC$, $\chi(\CC,a)$ is finite iff
its cardinality is smaller than $N$.
\end{claim}

Now, by Fact \ref{ext family}, it will sufficient to show that the extension property is
first-order for the family of complete types over $\emptyset$ that extends $\theta$ (for short we
will say that the extension property is first-order for $\theta$.) To see this, for any
$\chi(x,\bar y)\in L$ with $\chi(x,\bar y)\vdash \theta(x)$, consider the following $L_P$ formula
(by Claim \ref{H_claim} it is a formula):
$$S_\chi(x)=\forall\bar y\ [(P(\bar y)\wedge \exists^{<\infty}x\chi(x,\bar y))\rightarrow
\neg\chi(x,\bar y)].$$ For every finite set $\Delta$ of formulas of the form $\chi(x,\bar y)$ such
that $\chi(x,\bar y)\vdash \theta(x)$ ($x$ is the fixed variable of $\theta$, $\bar y$ any tuple of
variables), and consistent $\phi(x)\vdash\theta(x)$, let
$$\Psi_{\Delta,\phi}=\exists x[\phi(x)\wedge\bigwedge_{\chi\in\Delta}S_\chi(x)].$$

The following two claims shows that the extension property is first-order for $\theta$:

\begin{claim}
For any pair $\hat M=(M,P^M)$ of $T$ that satisfies the extension property for $\theta$, we have
$\hat M\models\Psi_{\Delta,\phi}$ for any finite set $\Delta$ (as above) and any $\phi(x)\vdash
\theta(x)$.
\end{claim}

\begin{claim}
Assume $\hat M=(M,P^M)$ is a pair of $T$ that is $\vert T\vert^+$-saturated and satisfies
$\Psi_{\Delta,\phi}$ for any finite set $\Delta$ (as above) and any $\phi(x)\vdash \theta(x)$. Then
$\hat M=(M,P^M)$ satisfies the extension property for $\theta$.
\end{claim}

\qedend

\section{Type-definability of the D-rank and more}

Recall the following lemma [S1, Lemma 8.4].

\begin{fact}\label{tilde-tau-lemma}
Assume the extension property is first-order in $T$. Let $\psi(x,z_1,...,z_m)$ be a Stone-open
relation over $\emptyset$ and let $\phi_j(x,y_j)\in L$ for $j=0,..,m$. Let $U$ be the following
invariant set. For all $d_1\in \CC$, $U(d_1)$  iff

$$\exists a\exists d_2...d_m[\psi(a,d_1,...d_m)\wedge\bigwedge_{j=0}^m (\phi_j(a,y_j) {\mbox{ forks\
over\ } d_1...d_j})].$$ Then $U$ is a $\tau^f$-open set over $\emptyset$. If each $\phi_j(x,y_j)$
is assumed to be low in $y_j$ and $\psi$ is assumed to be definable, then $U$ is a basic
$\tau^f_{\infty}$-open set.
\end{fact}

The following variation of the above fact will be useful to us.

\begin{lemma}\label{tilde-tau-lemma1}
Assume the extension property is first-order in $T$. Let $\psi(x,z_1,...,z_m)$ be a Stone-open
relation over $\emptyset$ and let $\phi_j(x,y_j)\in L$ for $j=1,..,m$. Let $U$ be the following
invariant set. For all $d_1\in \CC$, $U(d_1)$  iff

$$\exists a\exists d_2...d_m[\psi(a,d_1,...d_m)\wedge\bigwedge_{j=1}^m (\phi_j(a,y_j) {\mbox{ forks\
over\ } d_1...d_j})].$$

\noindent Then $U$ is a Stone-open set over $\emptyset$. Moreover, if we assume in addition that
each $\phi_j(x,y_j)$ is low in $y_j$ and $\psi(x,z_1,...,z_m)$ is definable, then $U$ is a
definable set over $\emptyset$.
\end{lemma}

\begin{remark}\em
The point in Lemma \ref{tilde-tau-lemma1} is that the index $j$ starts from 1 instead of 0 in Fact
\ref{tilde-tau-lemma}.
\end{remark}

\proofof {Lemma \ref{tilde-tau-lemma1}} First, the proof of the following claim is identical to the
proof of [S1, Subclaim 8.5].

\begin{claim}\label{Q generalized}
Let $\Gamma'$ be defined by $\Gamma'(d_1)$ iff $$\bigwedge_{\bar\eta=\{\eta_j\}_{j<m}\in L} \forall
d_2...d_m[(\bigwedge_{j=1}^{m-1} \eta_j(d_1...d_m,y_j)\ \mbox{forks\ over}\ d_1...d_j)\rightarrow
\forall a \Lambda_{\bar\eta}(a,d_1,...,d_m)].$$ where $\Lambda_{\bar\eta}$ is defined by
$$\Lambda_{\bar\eta}(a,d_1,...d_m)\ \ \mbox{iff}\ \ \psi(a,d_1,...d_m)\rightarrow \bigvee_{j=1}^m [\phi_j(a,y_j)\wedge
\neg\eta_j(d_1...d_m,y_j)\ \mbox{dnfo}\ d_1...d_m],$$ where $\eta_m$ denotes a contradiction. Then
$\Gamma'=\Gamma$ (Note that when $m=1$ the above empty conjunction from $j=1$ to $m-1$ is
interpreted as "True").
\end{claim}

\noindent Note that since the extension property is first-order in $T$, the relation
$\Lambda_{\bar\eta}^0$ defined by $\Lambda_{\bar\eta}^0(d_1,...d_m)\equiv\forall a
\Lambda_{\bar\eta}(a,d_1,...d_m)$ is type-definable by Fact \ref{foext}.

\noindent Thus we conclude that if $\Gamma$ is the complement of $U$ then for all $d_1$,
$\Gamma(d_1)$ iff

$$\bigwedge_{\bar\eta=\{\eta_j\}_{j<m}\in L}\forall d_2...d_m[\neg\Lambda_{\bar\eta}^0(d_1,...,d_m)\rightarrow
\bigvee_{j=1}^{m-1} \eta_j(d_1...d_m,y_j)\ \mbox{dnfo}\ d_1...d_j].$$

Thus, we conclude that if $m=1$, $\Gamma$ is type-definable (as the disjunction in the last formula
from $j=1$ to $m-1$ is interpreted as "False") and if $m>1$, we are done by the induction
hypothesis. The "Moreover" claim follows easily by Remark \ref{low type def} and compactness.

\begin{corollary}\label{type_def_d_rank}
Let $T$ be any simple theory in which that extension property is first-order. Let $\phi(x,y)\in L$
and let $k<\omega$. Then the set $$D_\phi^{\leq k}\equiv\{a \vert\ D(\phi(x,a))\leq k \}$$ is a
type-definable set (over $\emptyset$).
\end{corollary}

\proof Straightforward by \ref{tilde-tau-lemma1}.

\section{Continuity of the SU-rank}
In this section we work in $\CC$ unless stated otherwise and all tuples $a,b,c...$ are assumed to
be finite.

\begin{notation}
For a formula $\theta\in L$ and finite tuple of sorts $s$, let $$An^s(\theta)=\{a\in \CC^s\vert\
tp(a)\ \mbox{is\ analyzable\ in}\ \theta\}.$$ For finite tuples sorts $s_0,s_1$, let
$An^{s_0,s_1}(\theta)=An^{s_0}(\theta)\times An^{s_1}(\theta)$.
\end{notation}

\begin{fact}\label{An_open}
For every formula $\theta\in L$, and tuple of sorts $s$, $An^s(\theta)$ is a Stone-open set.
\end{fact}

\begin{lemma}\label{main_lemma}
Let $T$ be a simple theory in which the extension property is first-order. Assume $\theta(x)\in L$
and $SU(\theta)=1$. Let $r<\omega$ and let $s_0,s_1$ be finite tuple of sorts. Then the set
$$SU^{Aint}_{\leq r,s_0,s_1}(\theta)\equiv\{(a,a')\in \CC^{s_0}\times \CC^{s_1}\vert\ SU(a/a')\leq r\ \mbox{and\ } tp(a/a')\ \mbox{is\ almost\ internal\ in\ }
\theta\}$$ is a Stone-open set (over $\emptyset$).
\end{lemma}

\proof Assume $(a,a')\in SU^{Aint}_{\leq r,s_0,s_1}(\theta)$, i.e. $tp(a/a')$ is almost
$\theta$-internal and $SU(a/a')=k\leq r$. By Fact \ref{fact_a_int}, there exists $b\supseteq a'$
such that $b\backslash a'$ is a tuple of realizations of $tp(a)$, and there exists a formula
$\chi(x,\bar y,b)$ such that $\forall \bar y \exists^{<\infty} x \chi(x,\bar y,b)$ and such that
for all $\hat a\models Lstp(a/a')$, there is a tuple $\bar c=c_1,...,c_n$ of realizations of
$\theta$ such that $\models \chi(\hat a,\bar c,b)$. Let $s$ be the sequence of sorts of $b$. Let
$\lambda=(2^{\vert T\vert})^+$. Let $U$ be the subset of $\CC^{s_0}\times \CC^{s_1}$ defined by
$$U(\hat a,\hat a')\ \mbox{iff}\ \exists \{b_i\}_{i<\lambda}, b_i\in\CC^s,  b_i\supseteq \hat a'\ [\forall \tilde a\models tp(\hat a/\hat
a')\ \exists \bar c\subseteq\theta^\CC (\bigvee_{i<\lambda} \chi(\tilde a,\bar c,b_i))].$$

\noindent Now, for $l\leq n$, let $\delta_{l,n}$ be defined in the following way:
$\delta_{l,n}(\bar d,e)$ iff $\bar d=d_1,...d_n$, where each $d_i$ has the sort of a single
realization of $\theta$ and for some distinct $1\leq i_1,i_2,...,i_l\leq n$, $\bigwedge_{j=1}^l
d_{i_j}\not\in acl(e,d_{i_1}...d_{i_{j-1}})$. Note that if $\bar d=d_1,...,d_n$ and each $d_i$
realizes $\theta$, then for all $e$, $\delta_{l,n}(\bar d,e)$ iff $dim(\bar d/e)\geq l$ iff
$SU(\bar d/e)\geq l$.

\noindent In addition, let $\Theta$ be the subset of $\CC^{s_0}\times \CC^{s_1}$ defined by
$$\Theta(\hat a,\hat a')\ \mbox{iff}\ \forall b'\in \CC^s [b'\supseteq \hat a'\rightarrow
\bigwedge_{l=k+1}^n (\Phi_l(\hat a,b')\rightarrow \Psi_{l-k}(\hat a,b'))] $$

\noindent where, $\Phi_l$, for $l<\omega$, is the subset of $\CC^{s_0}\times \CC^s$ defined by

$$\Phi_l(\hat a,b') \mbox{iff}\ \exists\bar c=c_1,c_2,...c_n\in \theta^\CC\ [\chi(\hat
a,\bar c,b')\wedge \delta_{l,n}(\bar c,b')],$$ and $\Psi_l$, for $l<\omega$, is the subset of
$\CC^{s_0}\times \CC^s$ defined by

$$\Psi_l(\hat a,b') \mbox{iff}\ \exists\bar c=c_1,c_2,...c_n\in \theta^\CC\ [\chi(\hat
a,\bar c,b')\wedge \delta_{l,n}(\bar c,\hat ab')].$$

\begin{subclaim}
$(a,a')\in U\cap\Theta$.
\end{subclaim}

\proof Since the number of extensions of $tp(a/a')$ to $bdd(a')$ is small, we get the required
$b_i$-s in the definition of $U$ for ensuring that $U(a,a')$ holds. Now, to show that
$(a,a')\in\Theta$, assume $b'\in\CC^s$ and $b'\supseteq a'$ and $\Phi_l(a,b')$ holds for some
$k+1\leq l\leq n$. Then, there is $\bar c=c_1,c_2,...c_n\in \theta^\CC$ such that $\models
\chi(a,\bar c,b')$ and $SU(\bar c/b')\geq l$. Thus $SU(\bar c/b')=SU(\bar ca/b')\leq SU(a/b')\oplus
SU(\bar c/ab')$. Since we know that $SU(a/b')\leq k$ (as $SU(a/a')=k$ and $b'\supseteq a'$), we
conclude that $SU(\bar c/ab')\geq l-k$.

\begin{subclaim}
$U\cap\Theta\subseteq SU^{Aint}_{\leq r,s_0,s_1}(\theta)$.
\end{subclaim}

\proof Assume $(\hat a,\hat a')\in U\cap\Theta$. Since $(\hat a,\hat a')\in U$, by extension there
exists $b^*\in \CC^s$ such that $b^*\supseteq \hat a'$ and $\nonfork {b^*}{\hat a}{\hat a'}$ and
such that for some $\bar c=c_1,c_2,...c_n\in \theta^\CC$ we have $\models \chi(\hat a,\bar c,b^*)$.
We may assume that $SU(\bar c/b^*)\geq SU(\bar c'/b^*)$ for all tuples $\bar
c'=c'_1,c'_2,...c'_n\in \theta^\CC$ such that $\models \chi(\hat a,\bar c',b^*)$. Let $l=SU(\bar
c/b^*)$. Then $\models\Phi_l(\hat a,b^*)$. By $\models\Theta(\hat a,\hat a')$, we conclude that
$\models\Psi_{l-k}(\hat a,b^*)$. So, there exists a tuple $\bar c^*=c^*_1,c^*_2,...c^*_n\in
\theta^\CC$ such that $\models \chi(\hat a,\bar c^*,b^*)$ and $SU(\bar c^*/\hat ab^*)\geq l-k$.
Now, as $SU(\hat a\bar c^*/b^*)<\omega$,
$$SU(\bar c^*/b^*)=SU(\hat a\bar c^*/b^*)=SU(\hat a/b^*)+SU(\bar c^*/\hat ab^*).$$ By maximality of
$l$, $SU(\bar c^*/b^*)\leq l$; so it follows that $SU(\hat a/b^*)\leq k$. Since $\nonfork
{b^*}{\hat a}{\hat a'}$ and $b^*\supseteq \hat a'$, we conclude that $SU(\hat a/\hat a')\leq k$.

\begin{subclaim}
$U\cap\Theta$ is Stone-open.
\end{subclaim}

\proof It will be sufficient to prove that $\Psi_l$ is Stone-open for every $l\leq n$. Indeed, let
us define for every distinct $1\leq i_1,i_2,...,i_l\leq n$, a relation $\Psi_{i_1,i_2,...,i_l}$ by:
for every $(\hat a,b')\in\CC^{s_0}\times \CC^{s_1}$, we have $$\Psi_{i_1,i_2,...,i_l}(\hat
a,b')\mbox{\ iff\ }\exists\bar c=c_1,c_2,...c_n\in \theta^\CC\ [\chi(\hat a,\bar c,b')\wedge
\bigwedge_{j=1}^l (y=c_{i_j}\mbox{\ forks\ over\ }\hat ab'c_{i_1}...c_{i_{j-1}})].$$ Then, for
every $(\hat a,b')\in\CC^{s_0}\times \CC^{s_1}$, we have $\Psi_l(\hat a,b')$ iff there are distinct
$1\leq i_1,i_2,...,i_l\leq n$ such that $\Psi_{i_1,i_2,...,i_l}(\hat a,b')$. By Lemma
\ref{tilde-tau-lemma1}, $\Psi_{i_1,i_2,...,i_l}$ is definable for every distinct $1\leq
i_1,i_2,...,i_l\leq n$. To see this, choose in Lemma \ref{tilde-tau-lemma1}: $m=l$, $d_1=\hat ab'$,
$d_j$ remains the same for $j=2,...l$, $a=c_1c_2...c_n$, $\phi_j=(y=c_{i_j})$, and
$$\psi(a,d_1,...,d_m)=\chi(\hat a,\bar c,b')\wedge\bigwedge_{j=2}^l
(d_j=c_{i_{j-1}})\wedge\bigwedge_{i=1}^n \theta(c_i).$$

\begin{lemma}\label{lemma1}
Let $T$ be any simple theory. Assume $\theta(x)\in L$ and $SU(\theta)=1$. Let $r<\omega$ and let
$s_0,s_1$ be sorts. Then the set
$$SU^{An,Aint}_{\leq r,s_0,s_1}(\theta)\equiv\{(a,a')\in An^{s_0,s_1}(\theta) \vert\ SU(a/a')\leq r\ \mbox{and\ } tp(a/a')\ \mbox{is\ almost\ internal\ in\ }
\theta\}$$ is a Stone-open set (over $\emptyset$).
\end{lemma}

\proof Let $(a,a')\in SU^{An,Aint}_{\leq r,s_0,s_1}(\theta)$, i.e. $tp(a/a')$ is almost
$\theta$-internal, $tp(a')$ is analyzable in $\theta$ and $SU(a/a')\leq r$. Let
$(a_1,a_2,...a_{n-1},a_n)$ be an a-analysis in $\theta$ with $a_{n-1}=a',a_n=a$. Let
$\UU(x_1,...,x_n)$ be the Stone-open set over $\emptyset$ such that $(a'_1,...a'_n)\models\UU$ iff
$(a'_1,a'_2,...a'_n)$ is an a-analysis in $\theta$. Let $\phi(x_1,...,x_n)\in L$ be such that
$(a_1,a_2,...a_{n-1},a_n)\models\phi$ and $\phi^\CC\subseteq\UU$. Let $V_0=\theta$ and for $1\leq
i\leq n$ let $V_i$ be the projection of $\phi(x_1,...,x_n)$ on the $i-th$ coordinate. For $0\leq
i\leq n$ let $\hat s_i$ be the sort of $V_i$. By working in $\Ceq$ we may assume that for all
$0\leq i<j\leq n$, $\hat s_i\neq \hat s_j$. Let $\bar M=(V_0,V_1,...,V_n)$ be the structure whose
universe is the disjoint union of the $V_i$-s, where the interpretation of $\hat s_i$ in $\bar M$
is $V_i$ and $\bar M$ is equipped with the induced structure from $\CC^{eq}$, that is, the
$\emptyset$-definable subsets of $\bar M$ are precisely the $\emptyset$-definable sets of
$\CC^{eq}$ that are subsets of the cartesian products
of $V_{i}^{k_i}$ for some $i$-s and $k_i$-s. Clearly $\bar M$ is saturated.\\

\noindent In the following, if $p(x)$ is a partial type of $\bar M$, we will consider this type in
$\CC^{eq}$ by replacing $p(x)$ with $p'(x)=p(x)\wedge \nu(x)$ where $\nu(x)$ is the formula that
says "$x$ belongs to $\bar M$"

\begin{claim}\label {claim1}
\noindent 1) Dividing of partial types of $\bar M$ is absolute between $\bar M$ and $\CC^{eq}$ and
$\bar T=Th(\bar M)$ is simple (thus forking is is absolute between $\bar M$ and $\CC^{eq}$).

\noindent 2) Almost internality for types in $\bar M$ is absolute between $\CC^{eq}$ and $\bar M$.

\noindent 3) Every type of $\bar M$ is a-analyzable in $V_0=\theta$ in $\bar M$ (and in $\Ceq$).
Clearly, $\theta$ is weakly minimal in $\bar M$ too, thus $\bar T$ is supersimple.
\end{claim}

\proof 1) Clearly dividing for partial types of $\bar M$ is absolute between $\bar M$ and
$\CC^{eq}$. Thus, it is clear that every complete finitary type of $\bar M$ does not divide over a
subset of size $\leq\vert T\vert$ of its domain, so $\bar T$ is simple. To prove 2), assume that
$p\in S(A)$ and $q$ is a partial type over $A$ both of $\bar M$. By 1) if $p$ is almost internal in
$q$ in the sense of $\bar M$ then the same holds in the sense of $\CC^{eq}$. Now, assume $p$ is
almost internal in $q$ in the sense of $\Ceq$, where $p$ and $q$ are types of $\bar M$ ($p$ is
complete). By Fact \ref{fact_a_int}, in $\Ceq$ there exists $a\models p$ and a tuple $\bar a$ of
realizations of $p$ that is independent from $a$ such that $a\in acl(q^\CC,\bar a)$ (over the
corresponding parameters). Since $a,\bar a\in\bar M$, using 1), the same is true in $\bar M$. To
prove 3) we need to prove that for every $a,A\subseteq\bar M$, $tp(a/A)$ is a-analyzable in $V_0$
in $\bar M$. By Fact \ref{some internal}, we may assume that $a$ is a singleton (rather than a
tuple) and $A=\emptyset$. If $a\in V_0$ we are clearly done. Otherwise, $a\in V_i$ for some $1\leq
i\leq n$, so there are $a_1,...,a_{i-1},a_{i+1},...a_n\in\CC$ such that $\Ceq\models
\phi(a_1,...,a_{i-1},a,a_{i+1},...a_n)$. Then, in particular, in $\Ceq$, $(a_1,...,a_{i-1},a)$ is
an a-analysis in $\theta$. By 2), $(a_1,...,a_{i-1},a)$ is an a-analysis in $\theta$ in $\bar
M$.\qed

\begin{claim}\label{claim2}
For every $a,A$ in $\bar M$ we have $SU^{\bar M}(a/A)=SU(a/A)$, where $SU^{\bar M}$ is the
$SU$-rank in $\bar M$ and $SU$ is the usual $SU$-rank in $\Ceq$.
\end{claim}

\proof It will be sufficient to show that $SU^{\bar M}(a/A)\geq SU(a/A)$ for every
$a,A\subseteq\bar M$. We prove by induction on $\alpha$ that for all $a,A$, if $SU(a/A)\geq\alpha$
then $SU^{\bar M}(a/A)\geq\alpha$. We may clearly assume $\alpha=\beta+1$. We work in $\Ceq$. Let
$B\supseteq A$ be such that $\fork{a}{B}{A}$ and $SU(a/B)\geq\beta$. Let $I$ be Morley sequence of
$Lstp(a/B)$ and let $a'\models Lstp(a/B)$ be such that $\nonfork{a'}{I}{bdd(B)}$ (so, $a'\in\bar
M$). Let $e=Cb(a/bdd(B))$. We conclude that $\fork{a'}{AI}{A}$ and
$SU(a'/AI)=SU(a'/Ae)=SU(a'/B)=SU(a/B)$ (as $\nonfork{a'}{BI}{Ae}$ and $e\in dcl(I)$). So,
$SU(a'/AI)\geq\beta$. As $AI\subseteq\bar M$, by Claim \ref{claim1} 1), we get that $\bar M\models
\fork{a'}{AI}{A}$. Thus $SU^{\bar M}(a/A)=SU^{\bar M}(a'/A)\geq \alpha=\beta+1$. \qed

$\\$

\noindent Now, let $U_{\bar M}=(SU^{Aint}_{\leq r,\bar s_n,\bar s_{n-1}}(\theta))^{\bar M}$ i.e.
$U_{\bar M}$ is the set defined in Lemma \ref{main_lemma} in $\bar M$. By Corollary
\ref{epfo_su1_formula} and Claim \ref{claim1}, the extension property is first-order in $\bar
T=Th(\bar M)$. Thus, by Lemma \ref{main_lemma}, $U_{\bar M}$ is a Stone-open set over $\emptyset$
in $\bar M$ and in particular in $\CC^{eq}$. Since almost-internality is absolute between $\bar M$
and $\CC^{eq}$ by Claim \ref{claim1}, and $SU^{\bar M}=SU$ by Claim \ref{claim2}, $(a,a')\in
U_{\bar M}$. By Claims \ref{claim1}, \ref{claim2}, $U_{\bar M}\subseteq SU^{An,Aint}_{\leq
r,s_0,s_1}(\theta)$.\qed

\begin{notation}
For a formula $\theta=\theta(x,c)\in L(\CC)$ and sort $s$, let $\tilde {An}^s(\theta)=\{a\in
\CC^s\vert\ tp(a/c')\ \mbox{is\ analyzable\ in}\ \theta(x,c') \mbox{for\ all } \emptyset
\mbox{-conjugate}\ \theta(x,c') \mbox{of\ } \theta(x,c)\}.$ For sorts $s_0,s_1$, let $\tilde
{An}^{s_0,s_1}(\theta)=\tilde {An}^{s_0}(\theta)\times \tilde An^{s_1}(\theta)$.
\end{notation}

\begin{theorem}\label{cont_thm}
Let $T$ be a simple theory and work in $\Ceq$. Assume $\theta\in L(\CC)$, $SU(\theta)=1$. Let
$r<\omega$ and let $s_0,s_1$ be sorts. Then the set
$$SU^{\tilde{An}}_{\leq r,s_0,s_1}(\theta)\equiv\{(a,b)\in \tilde{An}^{s_0,s_1}(\theta) \vert\ SU(a/b)\leq r \}$$ is a Stone-open set (over $\emptyset$).
\end{theorem}

\proof First, assume that $\theta=\theta(x)\in L$. Note that $(a,b)\in SU^{An}_{\leq
r,s_0,s_1}(\theta)$ iff there are $a_0,a_1,...a_n=a\in dcl(ab)$ and $r_0,...,r_n<\omega$ such that
$r_0+...+r_n=r$, $tp(b)$ is analyzable in $\theta$ and $tp(a_i/a_{<i}b)$ is (almost)
$\theta$-internal and $SU(a_i/a_{<i}b)\leq r_i$ for all $i\leq n$. By Lemma \ref{lemma1} we are
done in case $\theta\in L$. Assume now $\theta=\theta(x,c^*)\in L(\CC)$. By working over $c^*$ and
what we have just proved, there is Stone-open set $U(x,y,c^*)$ over $c^*$ such that for all $(a,b)$
we have $U(a,b,c^*)$ iff $tp(a/c^*)$ and $tp(b/c^*)$ are analyzable in $\theta(x,c^*)$ and
$SU(a/bc^*)\leq r$. Let $p=tp(c^*)$ and let $V(x,y)=\forall c (c\models p\rightarrow U(x,y,c))$.
Then $V(x,y)$ is Stone-open over $\emptyset$ and $V^\CC=SU^{\tilde{An}}_{\leq r,s_0,s_1}(\theta)$
(whenever $(a,b)\models V$, there exists $c\models p$ independent from $ab$ over $\emptyset$ such
that $U(a,b,c)$; thus $SU(a/b)\leq r$.) \qed

\section{Definability of the D-rank}

\begin{theorem}
Let $T$ be a simple theory in which the extension property is first-order. Work in $\Ceq$. Assume
$\theta(x)\in L(\CC)$, $SU(\theta)=1$. Let $\phi(x,y)\in L$ be such that
$$\forall a,b\ [\phi(a,b)\Rightarrow a\in \tilde{An}^x(\theta), b\in \tilde{An}^y(\theta)].$$ Then the $D$-rank is
definable for $\phi(x,y)$, namely, for all $r<\omega$, the set $D_{\phi}^r=\{b\vert\
D(\phi(x,b))=r\}$ is definable.
\end{theorem}

\proof By Corollary  \ref{type_def_d_rank}, $D_{\phi}^{\leq r}=\{b\vert\ D(\phi(x,b))\leq r\}$ is
type-definable. Now, $D_{\phi}^{\leq r}=\{b\vert\ \forall a\ (\phi(a,b)\rightarrow\ SU(a/b)\leq
r)\}$. By Theorem \ref {cont_thm} and our assumption, $D_{\phi}^{\leq r}$ is Stone-open so we are
done.\\

We say that the $SU$-rank is \em uniformly continuous \em if for any given sorts $s_0,s_1$, the set
$\{(a,b)\in \CC^{s_0}\times \CC^{s_1}\vert\ SU(a/b)\leq r\}$ is Stone open for any $r$. We say that
\em the $D$-rank is definable \em if the $D$-rank is definable for $\phi(x,y)$ for any $\phi(x,y)$.

\begin{corollary}
Let $T$ be a supersimple unidimensional theory (e.g. $T$ is a countable hypersimple theory). Then
the $SU$-rank is uniformly continuous and the $D$-rank is definable. In particular, for every
finitary type $p\in S(A)$ we have $D(p)=SU(p)$.
\end{corollary}

Ziv Shami,  E-mail address: zivshami@gmail.com.

\end{document}